# A note on Waring's Problem

Li An-Ping

Beijing 100085, P.R. China
apli0001@sina.com

Abstract

In this paper, we will present a new iterative construction for the auxiliary equation of Waring's problem, which seems a little simpler than the one introduced in [4], and give the same upper bounds of $G(k)$ as the ones of the paper [4] and [8].

## 1. Introduction

Waring's problem is a well-known problem in number theory, the original statement is that, for each natural number $k$, there is a positive integer $g(k)$, such that each natural numbers may be represented a sum of at most $g(k)$ $k$th powers of natural numbers. Since the original problem was proved by Hilbert, and the values of $g(k)$ have been almost completely known, the modern version of the problem is to find the least integer $G(k)$ such that each sufficient large integer may be represented a sum of at most $G(k)$ $k$th powers of natural numbers. In about 1920's, Hardy, Ramanujan and Littlewood proposed so called circle method, which first appears hope in solving some problems of number theory, since then many mathematicians such as Vinogradov, Davenport, Hua, Chen, Vaughan, Wooley and others have made enduring and great efforts on exploiting Hardy-littlewood method, and have succeed in progress for $G(k)$. For the details is referred to see the Vaughan and Wooley's survey paper [6].

Suppose that $N$ is a sufficient large integer, for a given integer $k$, denoted by $P = \left[N^{1/k}\right]$. For $s$ integer sets $X_i, 1 \le i \le s, X_i \subseteq [0, P]$, denoted by $\mathfrak{A} = \left\{ \sum_{1\le i\le s} x_i^k \mid x_i \in X_i, 1 \le i \le s \right\}$. Let $\gamma(m)$ be the number of occurrences of number $m$ in $\mathfrak{A}$, that is, the characteristic function of set $\mathfrak{A}$. It is clear that $\prod_{1\le i\le s} |X_i| = \sum_m \gamma(m)$, hence

$$\left( \prod_{1\le i\le s} |X_i| \right)^2 = \left( \sum_m \gamma(m) \right)^2 \le \sum_{\gamma(m)>0} 1 \sum_m \gamma(m)^2$$

So,

$$\sum_{\gamma(m)>0} 1 \ge \left( \prod_{1\le i\le s} |X_i| \right)^2 \Big/ \sum_m \gamma(m)^2 .$$

The formula above indicates that $\sum_m \gamma(m)^2$ is more less, the size $|\mathfrak{A}|$ is more large when $\prod_{1\le i\le s} |X_i|$ given. In Hardy-littlewood method, an important work is to find as little as possible $s$ and sets $X_i, 1 \le i \le s$, such that $\sum_m \gamma(m)^2$ as small as possible and that set $\mathfrak{A}$ is dense in $[1, N]$ with some sense. Let $X = (X_1, X_2, \cdots, X_s)$, $S_s(P, X) = \sum_m \gamma(m)^2$, it is clear that

$S_s(P,X)$ is the number of solutions of equation

$$x_1^k+\cdots+x_s^k=y_1^k+\cdots+y_s^k,\quad x_i,y_i\in X_i, 1\le i\le s.$$

The equation above is called auxiliary equation of Waring problem.

Earlier work such as Vinogradov [7] and Davenport [1,2], the selection of sets $X_i, 1\le i\le s,$ are different from each other. Vaughan [4] creatively introduced a iterative method with the sets $X_i$ are same each other, $X_i=\mathscr{A}(P,R), 1\le i\le s,$ where $\mathscr{A}(P,R)$ is a set of so called "smooth" integers. In their papers, $S_s(P,X)$ is written as $S_s(P,R)$. It is not difficult to imagine that there is a merit for such homogeneous construction is that it may be represented as a simple integral form, which is favorable to be varied by Hölder inequality, this feature is reflected very well in the paper [4].

In this paper, we will present a new construction for domain set $X$ of the auxiliary equation, which is seemingly a little simpler than the one of "smooth" numbers in [4] and [8], but most arguments and results are almost same as ones in [4] and [8]. With the new homogeneous construction, we also obtain

**Theorem 1** (Vaughan 1989 [4] ). For sufficient large $k$,

$$G(k)\le 2k\left(\log(k\log k)+1+\log 2+O\left(\frac{\log\log k}{\log k}\right)\right). \tag{1.1}$$

**Theorem 2** (Wooley 1991 [8] ). For sufficient large $k$,

$$G(k)\le k\left(\log(k\log k)+O(1)\right). \tag{1.2}$$

**2. The Proof of Theorem 1.**

$P$ is sufficient great, $\theta$ is a constant, with that $\theta>1/k$. Let $\tilde{P}=P^{1+\theta}$, $\mathscr{P}$ is a set of prime numbers $p$ in interval $[P^\theta/2,P^\theta]$, write $|\mathscr{P}|=Z$, then we know $Z\doteq P^\theta/2\log(P^\theta)$. Define recursively

$$\mathscr{C}(\tilde{P})=\left\{x\cdot p \mid x\in\mathscr{C}(P), p\in\mathscr{P}\right\}. \tag{2.1}$$

In this paper, we simply write $S_s(P)=S_s(P,X)$ as $X=\mathscr{C}(P)$.

**Lemma 1**.

$$S_s(\tilde{P}) \ll Z^s S_s(P) + Z^{2s} P S_{s-1}(P). \tag{2.2}$$

Proof. As usual, write $e(x) = e^{2\pi i x}$, let

$$f(\alpha) = \sum_{x \in \mathscr{C}(P)} e\left(x^k \alpha\right), \ f(\alpha, p) = \sum_{x \in \mathscr{C}(P)} e\left(p^k x^k \alpha\right), \ \tilde{f}(\alpha) = \sum_{y \in \mathscr{C}(\tilde{P})} e\left(y^k \alpha\right).$$

Then clearly,

$$\tilde{f}(\alpha) = \sum_{p \in \mathscr{P}} \sum_{x \in \mathscr{C}(P)} e\left(p^k x^k \alpha\right) = \sum_{p \in \mathscr{P}} f(\alpha, p).$$

Applying Hölder's inequality, it has

$$S_s(\tilde{P}) = \int_0^1 \left| \tilde{f}(\alpha)^{2s} \right| d\alpha = \int_0^1 \left| \tilde{f}(\alpha)^s \tilde{f}(\alpha)^{-s} \right| d\alpha$$

$$= \int_0^1 \left( \sum_{p \in \mathscr{P}} \sum_{q \in \mathscr{P}} f(\alpha, p) f(-\alpha, q) \right)^s d\alpha$$

$$= \int_0^1 \left( \sum_{p=q} + \sum_{p \neq q} \right)^s d\alpha$$

$$\ll \int_0^1 \left( Z \cdot |f(\alpha)|^2 \right)^s d\alpha + \int_0^1 ( \sum_{p,q \in \mathscr{P},\ p \neq q} |f(\alpha, p) f(-\alpha, q)| )^s d\alpha$$

$$\ll Z^s S_s(P) + Z^{2(s-1)} \sum_{p,q \in \mathscr{P},\ p \neq q} \int_0^1 |f(\alpha, p) f(\alpha, q)|^s d\alpha.$$

Moreover, let $\Lambda(\alpha, p, q) = f(\alpha, p)^{s-1} f(\alpha, q)$, then by Cauchy inequality, it has

$$\int_0^1 \left| f(\alpha, p)^s f(\alpha, q)^s \right| d\alpha = \int_0^1 |\Lambda(\alpha, p, q)| |\Lambda(\alpha, q, p)| d\alpha$$

$$\leq \left( \int_0^1 |\Lambda(\alpha, p, q)|^2 d\alpha \right)^{1/2} \left( \int_0^1 |\Lambda(\alpha, q, p)|^2 d\alpha \right)^{1/2}$$

And

$$\sum_{p,q \in \mathscr{P}, p \neq q} \int_0^1 \left| f(\alpha, p)^s f(\alpha, q)^s \right| d\alpha \leq \sum_{p,q \in \mathscr{P}, p \neq q} \left( \int_0^1 |\Lambda(\alpha, p, q)|^2 d\alpha \right)^{1/2} \left( \int_0^1 |\Lambda(\alpha, q, p)|^2 d\alpha \right)^{1/2}$$

$$\leq \sum_{p,q \in \mathscr{P}, p \neq q} \int_0^1 |\Lambda(\alpha, p, q)|^2 d\alpha$$

Denoted by $T_{p,q}$ the inner integral, which is the number of solutions of equation

$$p^k \left( (x_1^k + x_2^k + \cdots + x_{s-1}^k) - (y_1^k + y_2^k + \cdots + y_{s-1}^k) \right) = q^k \left( y^k - x^k \right) \tag{2.3}$$

With $x, y \in \mathscr{C}(P),\ x_i, y_i \in \mathscr{C}(P), 1 \leq i \leq s-1,\ \ p, q \in \mathscr{P},\ p \neq q.$

Hence, it has that $y^k \equiv x^k \mod p^k$. By the definition, we can know that, for $x \in \mathscr{C}(P)$, it has $(x, p) = 1$, since the prime factors of element $x$ less than $p$. Moreover, we know that for an integer $m,\ p^k \nmid m$, there are at most finite solutions for the congruence $x^k \equiv m \mod p^k$. So, we may divide residue system $\mathbb{Z}_{p^k}$ into finite classes $\mathbb{Z}_{p^k} = \bigcup_{1\le i\le l} A_i,\ (l \le k)$, such that there is at most one solution in each class $A_i$ for the congruence. This means, in one class $A_i$, $y^k \equiv x^k \mod p^k$ implies $y \equiv x \mod p^k$, and so here $x = y$, for $\theta > 1/k$. In more detail, let

$$\mathscr{C}_i(P) = \{\, x \mid x \in \mathscr{C}(P), \overline{x} \in A_i \,\},\ \ f_i(\alpha, q) = \sum_{x\in\mathscr{C}_i(P)} e\left(q^k x^k \alpha\right),\ 1 \le i \le l,$$

where $\overline{x}$ represents the least non-negative integer with that $x \equiv \overline{x} \mod p^k$.

Then, $f(\alpha, q) = \sum_{1\le i\le l} f_i(\alpha, q)$, and

$$T_{p,q} = \int_0^1 |\Lambda(\alpha, p, q)|^2\, d\alpha = \int_0^1 |f(\alpha, p)|^{2(s-1)} |f(\alpha, q)|^2\, d\alpha$$

$$= \int_0^1 |f(\alpha, p)|^{2(s-1)} \left|\sum_{1\le i\le l} f_i(\alpha, q)\right|^2 d\alpha$$

$$\le \int_0^1 |f(\alpha, p)|^{2(s-1)}\, l \cdot \sum_{1\le i\le l} |f_i(\alpha, q)|^2\, d\alpha$$

$$= l \cdot \sum_{1\le i\le l} \int_0^1 |f(\alpha, p)|^{2(s-1)} |f_i(\alpha, q)|^2\, d\alpha$$

Clearly, the integral in the summation is the number of solutions of equation

$$p^k \left((x_1^k + x_2^k + \cdots + x_{s-1}^k) - (y_1^k + y_2^k + \cdots + y_{s-1}^k)\right) = q^k \left(y^k - x^k\right) \tag{2.4}$$

With $x_j, y_j \in \mathscr{C}(P), 1 \le j \le s-1,\ x, y \in \mathscr{C}_i(P)$.

Now for the equation (2.4), it has $y^k \equiv x^k \mod p^k \Rightarrow y \equiv x \mod p^k \Rightarrow y = x$, for $\theta > 1/k$.

Therefore,

$$T_{p,q} \ll P S_{s-1}(P),$$

and

$$\sum_{p,q\in\mathscr{P},\ p\ne q} \int_0^1 |f(\alpha, p) f(\alpha, q)|^s\, d\alpha \ll Z^2 P S_{s-1}(P).$$

And then,

$$S_s(\tilde{P}) \ll Z^s S_s(P) + Z^{2s} P S_{s-1}(P). \qquad \square$$

**Corollary 1**. Suppose that $S_s(P) = P^{\lambda_s}$, then for $s \ge 2$,

$$\lambda_s \le (2s-k)+(k-2)\left(\frac{k}{k+1}\right)^{s-2}. \qquad (2.5)$$

Proof. By Lemma 1, it has
$\tilde{P}^{\lambda_s} \ll Z^s P^{\lambda_s} + Z^{2s} P^{\lambda_{s-1}+1}$.

As $\lambda_s \ge s$, and $Z \doteq P^{\theta}/(2\theta \log P)$, so the first term of the right-side of the inequality above is sub-term may be omitted, and hence

$P^{(1+\theta)\lambda_s} \ll P^{\lambda_{s-1}+1+2s\theta}$.

i.e.

$$\lambda_s \le \frac{\lambda_{s-1}}{(1+\theta)} + \frac{(1+2s\theta)}{(1+\theta)}. \qquad (2.6)$$

From the recursive inequality, and the known result $\lambda_2 = 2$, with a simple calculate, it follows

$$\lambda_s \le 2s - \frac{1}{\theta} + \left(\frac{1}{\theta} - 2\right)\left(\frac{1}{1+\theta}\right)^{s-2}$$

Let $\theta \to 1/k$, (5) is followed. □

Let $\tau = 2kP^{k-1}$, $\mathfrak{I} = [\tau^{-1}, 1+\tau^{-1}]$. For integers $a, q$, denotes

$\mathfrak{M}(q,a) = \left\{\alpha \,\middle|\, |\alpha - a/q| \le 1/q\tau\right\}$,

And for positive number $W$ $(\le P)$, denotes

$\mathfrak{N}(q,a) = \left\{\alpha \,\middle|\, |\alpha - a/q| \le W/(q\tau P)\right\}$.

$\mathfrak{M}$ is the union of the $\mathfrak{M}(q,a)$ with $1 \le a \le q \le P$, $(a,q)=1$. $\mathfrak{m} = \mathfrak{I} \setminus \mathfrak{M}$.

$\mathfrak{N}$ is the union of the $\mathfrak{N}(q,a)$ with $1 \le a \le q \le W$, $(a,q)=1$. $\mathfrak{n} = \mathfrak{I} \setminus \mathfrak{N}$.

Let $f(\alpha) = \sum_{x \le P} e(\alpha x^k)$, $g(\alpha) = \sum_{x \in \mathscr{C}(P)} e(\alpha x^k)$.

The proof of Theorem 1 will require the following results.

**Lemma 2.** For $\alpha \in \mathfrak{m}$,

$$f(\alpha) \ll P^{1-1/2^{k-1}}. \qquad (2.7)$$

**Lemma 3.** Suppose that $k \ge 3$, $s \ge k+2$, then

$$\int_{\mathfrak{M}} \left| f(\alpha)^s \right| d\alpha \ll P^{s-k}, \tag{2.8}$$

$$\int_{\mathfrak{M}\setminus\mathfrak{n}} \left| f(\alpha)^s \right| d\alpha \ll W^{\varepsilon-1/k} P^{s-k}. \tag{2.9}$$

**Lemma 4.** Let

$$V(\alpha) = \sum_{x/2<p\le X} \sum_{y\le Y} b_y e(yp^k\alpha),$$

where $b_y$ are arbitrary complex numbers. Suppose that $\alpha = a/q+\beta,$ with $|\beta| \le \frac{1}{2} q^{-1} X^{-k}$, $q \le 2X^k, (a,q)=1$, that $Y \gg X^k$, and that when $q \le X$ one has $|\beta| \gg q^{-1} X^{1-k} Y^{-1}$. Then

$$V(\alpha) \ll \left( XY^{1+\varepsilon} \sum_{y\le Y} |b_y|^2 \right)^{1/2}. \tag{2.10}$$

(2.7) is from Weyl' inequality, Lemmas 3 and 4 are just Lemma 5.4 of [4] and Lemma 5.1 of [3] respectively, the proofs are referred to see [3], [4], or [7].

**Lemma 5.** Let $X = P^{1/2}$, $h(\alpha) = \sum_{(X/2)\le p\le X} \sum_{x\in\mathscr{C}(X)} e(\alpha p^k x^k)$. Then

$$|h(\alpha)| \ll P^{1-\hat{\sigma}}, \quad \text{for } \alpha \in \mathfrak{n}. \tag{2.11}$$

where $\hat{\sigma} = \dfrac{\log(1+1/k)}{4(1+\lambda)}$, $\lambda = \log k + \log\log k + O\left(\dfrac{\log\log k}{\log k}\right)$.

Proof. By Hölder inequality, it has

$$|h(\alpha)|^{2s} \le X^{2s-1} \sum_{(X/2)\le p\le X} \left| \sum_{y\in Y} b_y e(\alpha p^k y) \right|^2$$

Where $Y = sX^k$, and $b_y$ is the number of solutions of

$$x_1^k + x_2^k + \cdots + x_s^k = y, \quad x_i \in \mathscr{C}(X), 1\le i\le s.$$

When $\alpha \in \mathfrak{n}$, by lemma 4 and (2.5), it has

$$|h(\alpha)|^{2s} \ll X^{2s-1+k+\lambda_s} \le X^{4s-1+(k-2)(k/(k+1))^{s-2}}$$

That is,

$$|h(\alpha)| \ll P^{(4s-1+(k-2)(k/(k+1))^{s-2})/4s} = P^{1-(1-(k-2)(k/(k+1))^{s-2})/4s} = P^{1-\sigma}$$

Where $\sigma = \sigma(k,s) = (1-(k-2)(k/(k+1))^{s-2})/4s$

From calculus, we know that when $s = \lambda/\log(1+1/k)$, $\sigma(k,s)$ arrives the maximum value

$$\hat{\sigma} = \frac{\log(1+1/k)}{4(1+\lambda)},$$

Where $\lambda$ is the root of equation $\left(1+\lambda\right)\beta = e^{\lambda}$, $\beta = (k-2)(k+1)^2/k^2$. It is easy to know

$$\lambda = \log k + \log\log k + O\left(\frac{\log\log k}{\log k}\right). \qquad \square$$

Suppose $u$ and $t$ are two natural numbers, with that $u \ge k+1$, $2t\hat{\sigma} \ge (k-2)\left(k/(k+1)\right)^{u-2}$,

define

$$L = \int_0^1 \left| f(\alpha)^2 g(\alpha)^{2u} h(\alpha)^{2t} \right| d\alpha, \quad I = \int_{\mathfrak{n}} \left| f(\alpha) g(\alpha)^{2u} h(\alpha)^{2t} \right| d\alpha.$$

**Lemma 6.**

$$L \ll P^{2+2t+2u-k}, \quad I \ll P^{1+2t+2u-k} W^{-\delta}. \tag{2.12}$$

Proof.

$$L = \int_{\mathfrak{m}} \left| f(\alpha)^2 g(\alpha)^{2u} h(\alpha)^{2t} \right| d\alpha + \int_{\mathfrak{M}} \left| f(\alpha)^2 g(\alpha)^{2u} h(\alpha)^{2t} \right| d\alpha.$$

By Hölder inequality and Lemmas 1, 5, it has

$$\int_{\mathfrak{m}} \left| f(\alpha)^2 g(\alpha)^{2u} h(\alpha)^{2t} \right| d\alpha \ll P^{2+2t(1-\hat{\sigma})} \int_{\mathfrak{m}} \left| g(\alpha)^{2u} \right| d\alpha$$

$$\ll P^{2+2t(1-\hat{\sigma})} P^{2u-k+(k-2)(k/(k+1))^{u-2}} \ll P^{2+2t+2u-k}.$$

$$\int_{\mathfrak{M}} \left| f(\alpha)^2 g(\alpha)^{2u} h(\alpha)^{2t} \right| d\alpha \le \left( \int_{\mathfrak{M}} \left| f(\alpha)^{k+2} \right| d\alpha \right)^{2/(k+2)} \left( \int_0^1 \left| g(\alpha)^{2u} h(\alpha)^{2t} \right|^{(k+2)/k} d\alpha \right)^{k/(k+2)}$$

$$\int_0^1 \left| g(\alpha)^{2u} h(\alpha)^{2t} \right|^{(k+2)/k} d\alpha = \int_0^1 \left| g(\alpha)^{2u+4u/k} h(\alpha)^{2t+4t/k} \right| d\alpha$$

$$= \int_0^1 \left| g(\alpha)^{2+2u} h(\alpha)^{2t} \right| \left| g(\alpha)^{4u/k-2} h(\alpha)^{4t/k} \right| d\alpha$$

$$\le P^{4u/k-2+4t/k} \int_0^1 \left| g(\alpha)^{2u+2} h(\alpha)^{2t} \right| d\alpha$$

$$\le P^{4u/k-2+4t/k} \int_0^1 \left| g(\alpha)^{2u+2} h(\alpha)^{2t} \right| d\alpha$$

Considering underlying equation, we can know

$$\int_0^1 \left| g(\alpha)^{2u+2} h(\alpha)^{2t} \right| d\alpha \le \int_0^1 \left| f(\alpha)^2 g(\alpha)^{2u} h(\alpha)^{2t} \right| d\alpha = L.$$

In addition, by (2.8), it has $\int_{\mathfrak{M}} \left| f(\alpha)^{k+2} \right| d\alpha \ll P^2$. So

$$L \ll P^{2+2t+2u-k} + P^{2\times 2/(k+2)} \cdot \left( L P^{(4u-2k+4t)/k} \right)^{k/(k+2)}$$

$$\ll P^{2+2t+2u-k} + P^{2\times 2/(k+2)} \cdot L^{k/(k+2)} \left(P^{(4u-2k+4t)/(k+2)}\right).$$

And,

$$L \ll P^{(2+2t+2u-k)}.$$

Similarly,

$$I = \int_{\mathfrak{n}} \left| f(\alpha) g(\alpha)^{2u} h(\alpha)^{2t} \right| d\alpha$$

$$= \int_{\mathfrak{m}} \left| f(\alpha) g(\alpha)^{2u} h(\alpha)^{2t} \right| d\alpha + \int_{\mathfrak{M}\setminus\mathfrak{N}} \left| f(\alpha) g(\alpha)^{2u} h(\alpha)^{2t} \right| d\alpha$$

$$\ll P^{1+2u+2t-k-\delta} + \left( \int_{(\mathfrak{M}\setminus\mathfrak{N})} \left| f(\alpha) \right|^{k+2} d\alpha \right)^{1/(k+2)} \left( \int_{(\mathfrak{M}\setminus\mathfrak{N})} \left| g(\alpha)^{2u} h(\alpha)^{2t} \right|^{(k+2)/(k+1)} d\alpha \right)^{(k+1)/(k+2)}$$

$$\ll P^{1+2u+2t-k-\delta} + \left( \int_{(\mathfrak{M}\setminus\mathfrak{N})} \left| f(\alpha) \right|^{k+2} d\alpha \right)^{1/(k+2)} \left( \int_{(\mathfrak{M}\setminus\mathfrak{N})} \left| g(\alpha)^{2u} h(\alpha)^{2t} \right|^{(k+2)/(k+1)} d\alpha \right)^{(k+1)/(k+2)}$$

$$\int_{(\mathfrak{M}\setminus\mathfrak{N})} \left| g(\alpha)^{2u} h(\alpha)^{2t} \right|^{(k+2)/(k+1)} d\alpha \le \int_0^1 \left| g(\alpha)^{2u} h(\alpha)^{2t} \right|^{(k+2)/(k+1)} d\alpha$$

$$\le LP^{2u/(k+1)+2t/(k+1)-2}$$

$$I \ll P^{1+2u+2t-k-\delta} + \left(P^2 W^{-\delta}\right)^{1/(k+2)} \left(LP^{2u/(k+1)+2t/(k+1)-2}\right)^{(k+1)/(k+2)}$$

It follows

$$I \ll P^{1+2t+2u-k} W^{-\delta}. \qquad \square$$

On the other hand, by elementary number theory, it is known the distribution of prime numbers in a larger interval is asymptotically equal $\bmod q$. Explicitly, denoted by $\pi(x)$ and $\pi(c,d)$ the numbers of prime numbers in the intervals $[0,x]$ and $[c,d]$ respectively, and for $0 < a < q \le P$, $(a,q)=1$, $\Delta = [c,d]$, let

$$\pi(x;q,a) = \sum_{\substack{p\equiv a \bmod q \\ p\le x}} 1, \quad \pi(\Delta;q,a) = \sum_{\substack{p\equiv a \bmod q \\ c\le p\le d}} 1, \quad \zeta(x;q,a) = \sum_{\substack{t\equiv a \bmod q \\ t\in\mathscr{C}(x)}} 1.$$

Then

$$\pi(x;q,a) = (1+o(1))\frac{\pi(x)}{\varphi(q)}, \qquad \pi(\Delta;q,a) = (1+o(1))\frac{\pi(c,d)}{\varphi(q)}.$$

From formulas above, it is easy to infer that

$$\zeta(x;q,a) = (1+o(1))\frac{|\mathscr{C}(x)|}{\varphi(q)}. \tag{2.13}$$

This indicates that the contribution from $\mathfrak{N}$ of $g(\alpha)$ and $h(\alpha)$ may be treated in the usual way as for the sequence of integers $\bmod q$ but with a proportion of $|\mathscr{C}(P)|/P$. Besides, by a simply calculate, it is easy to know

$$|\mathscr{C}(P)| \doteq \frac{P}{\log(P)^{(\eta+1)/2}} \cdot \left(\frac{k+1}{2}\right)^{\eta}, \tag{2.14}$$

where $\eta = k \log\log P$. So, $|\mathscr{C}(P)| > P^{1-\varepsilon}$. Hence, we have

$$\int_{\mathfrak{N}} f(\alpha) g(\alpha)^{2u} h(\alpha)^{2t} e(-N\alpha) d\alpha \gg P^{1+2u+2t-k}.$$

and

$$R(n) = \int_0^1 f(\alpha) g(\alpha)^{2u} h(\alpha)^{2t} e(-N\alpha) d\alpha \gg P^{1+2u+2t-k}.$$

Take $v = u-2, t = 1 + \left[\frac{k-2}{2\hat{\sigma}}\left(\frac{k}{k+1}\right)^{v}\right]$, it follows

$$G(k) \le 7 + 2v + 2\left[\frac{k-2}{2\hat{\sigma}}\left(\frac{k}{k+1}\right)^{v}\right].$$

From calculus, we can know that the right-side of the above arrives the maximum as $v$ tends to $\log(\mu(k-2)/2\hat{\sigma})/\mu, \quad \mu = \log((k+1)/k)$. Hence,

$$G(k) \le 7 + 2\log(\mu(k-2)/2\hat{\sigma})/\mu + 2\left[\frac{1}{\mu}\right]$$
$$\le 2k\left(\log(k\log k) + 1 + \log 2 + O\left(\frac{\log\log k}{\log k}\right)\right).$$

**3. The Proof of Theorem 2.**

Denoted by $\mathscr{P}[a,b]$ the set of prime numbers in the interval $[a,b]$. Let $\theta_1, \cdots, \theta_k$ be $k$ real numbers, $0 \le \theta_i \le 1/k, 1 \le i \le k,$ which will be determined later. Let $Z_i = P^{\theta_i}$, $P_{i+1} = P_i / Z_{i+1}$, $0 \le i < k$, and $H_i = P / Z_i^k$, $\mathscr{P}_i = \mathscr{P}[Z_i/2, Z_i]$, $1 \le i \le k$. Recursively define

$$\mathscr{C}(P_i) = \{x \cdot p \mid x \in \mathscr{C}(P_{i+1}), p \in \mathscr{P}_{i+1}, (p,x) = 1\}, \qquad i = 0, 1, \cdots k. \tag{3.1}$$

Simply write $\mathscr{C}_i = \mathscr{C}(P_i), \quad i = 0, 1, \cdots, k.$

In the following, it will be used the notation of difference of a function: As usual, for an integer coefficient of polynomial $\phi(x)$, recursively define the forward differences

$$\Delta_1(\phi(x), t) = \phi(x+t) - \phi(x),$$

$$\Delta_{i+1}(\phi(x), h_1, h_2, \cdots, h_i, h_{i+1}) = \Delta_1\left(\Delta_i(\phi(x), h_1, h_2, \cdots, h_i), h_{i+1}\right). \quad i = 1, 2, \cdots.$$

Suppose that $t = h \cdot m$, $m$ is a constant, then we know that $m \mid (\phi(x+t) - \phi(x))$, in this case we define modified differences

$$\Delta_1^*(\phi(x), h; m) = m^{-1}\left(\phi(x+hm) - \phi(x)\right).$$

$$\Delta_{i+1}^*(\phi(x), h_1, \cdots, h_i, h_{i+1}; m_1, \cdots, m_i, m_{i+1}) = \Delta_1^*\left(\Delta_i^*(\phi(x), h_1, \cdots, h_i; m_1, \cdots, m_i), h_{i+1}; m_{i+1}\right).$$

Simply write,

$$\Psi_i = \Psi_i(x, h_1, \cdots, h_i; p_1^k, \cdots, p_i^k) = \Delta_i^*(x^k, h_1, \cdots, h_i; p_1^k, \cdots, p_i^k), \qquad i = 1, 2, \cdots.$$

And define

$$F_i(\alpha, q) = \sum_{1 \le h_1 \le H_1} \cdots \sum_{1 \le h_i \le H_i} \sum_{p_1 \in \mathscr{P}_1} \cdots \sum_{p_i \in \mathscr{P}_i} \sum_x e(q^k \cdot \Psi_i \alpha), \qquad 1 \le i \le k. \tag{3.2}$$

$$f_i(\alpha) = \sum_{x \in \mathscr{C}_i} e(x^k \alpha), \quad f_i(\alpha, p) = \sum_{x \in \mathscr{C}_i, (x,p)=1} e(p^k x^k \alpha), \quad p \in \mathscr{P}_i, \ 1 \le i \le k. \tag{3.3}$$

Let

$$\begin{aligned} J_0 &= \sum_p T_{p,q} = \int_0^1 \left| f(\alpha, p)\right|^{2(s-1)} \left| f(\alpha, q)\right|^2 d\alpha \\ J_i &= \int_0^1 \left| f_i(\alpha)\right|^{2(s-1)} F_i(\alpha, q) d\alpha, \quad 1 \le i < k. \end{aligned} \tag{3.4}$$

Then, it has

**Lemma 7.**

$$\begin{aligned} J_0 &\le ZPS_{s-1}(P) + J_1, \quad J_i \ll U_i + V_i, \\ U_i &= S_{s-1}(P_i)^{1/2} Z_{i+1}^{(2s-3)/2} \left(P(\tilde{H}_i \tilde{Z}_i)^2 Z_{i+1} S_{s-1}(P_{i+1})\right)^{1/2} \\ V_i &= S_{s-1}(P_i)^{1/2} Z_{i+1}^{(2s-3)/2} (\tilde{H}_i \tilde{Z}_i J_{i+1})^{1/2}, \qquad 1 \le i < k, \end{aligned} \tag{3.5}$$

where $\tilde{H}_i = \prod_{j \le i} H_j$, $\tilde{Z}_i = \prod_{j \le i} Z_j$.

Proof. By Cauchy inequality and Hölder inequality,

$$J_i = \int_0^1 \left| f_i(\alpha)\right|^{2(s-1)} F_i(\alpha, q) d\alpha \le \left(\int_0^1 \left| f_i(\alpha)\right|^{2(s-1)} d\alpha\right)^{1/2} \left(\int_0^1 \left| f_i(\alpha)\right|^{2(s-1)} F_i(\alpha, q)^2 d\alpha\right)^{1/2}$$

$$\le \left(S_{s-1}(P_i)\right)^{1/2} \left(\int_0^1 \left| \sum_{p \in \mathscr{P}_{i+1}} f_{i+1}(\alpha, p) \right|^{2(s-1)} F_i(\alpha, q)^2 d\alpha\right)^{1/2}$$

$$\le \left(S_{s-1}(P_i)\right)^{1/2} \left(Z_{i+1}^{2s-3} \int_0^1 \sum_{p \in \mathscr{P}_{i+1}} \left| f_{i+1}(\alpha, p)\right|^{2(s-1)} F_i(\alpha, q)^2 d\alpha\right)^{1/2}$$

$$\le \left(S_{s-1}(P_i)\right)^{1/2} Z_{i+1}^{(2s-3)/2}$$

$$\times \left(\int_0^1 \sum_{p \in \mathscr{P}_{i+1}} \left| f_{i+1}(\alpha, p)\right|^{2(s-1)} P(\tilde{H}_i \tilde{Z}_i)^2 d\alpha + \left(\tilde{H}_i \tilde{Z}_i \int_0^1 \left| f_{i+1}(\alpha)\right|^{2(s-1)} F_{i+1}(\alpha, q)) d\alpha\right)\right)^{1/2}$$

$$\ll S_{s-1}(P_i)^{1/2} Z_{i+1}^{(2s-3)/2}$$
$$\times\left( \left( P(\tilde{H}_i \tilde{Z}_i)^2 Z_{i+1} S_{s-1}(P_{i+1}) \right)^{1/2} + \left( \tilde{H}_i \tilde{Z}_i \int_0^1 | f_{i+1}(\alpha)|^{2(s-1)} F_{i+1}(\alpha, q))d\alpha \right)^{1/2} \right)$$
$$\ll S_{s-1}(P_i)^{1/2} Z_{i+1}^{(2s-3)/2} \left( \left( P(\tilde{H}_i \tilde{Z}_i)^2 Z_{i+1} S_{s-1}(P_{i+1}) \right)^{1/2} + (\tilde{H}_i \tilde{Z}_i J_{i+1})^{1/2} \right)$$

□

It should be mentioned that in the proof of Lemma 7, we save the investigation for possible singular cases, that is, the cases $p \mid \Psi_i'(x)$. These cases may be treated as that divide the summation of $\Psi_i(x)$ into singular ones and normal part by Hölder inequality, the difference for the singular ones is delayed to the $i+1$ round, it will be showed the singular ones are secondary, and the number of singular ones is finite for a polynomial $\Psi_i$, hence it will be no effect for the result eventually.

In order that the differences can be carried out successively and efficiently to order $k$, it should be

$$U_i = V_i, \quad 0 \le i < k. \tag{3.6}$$

At first, we consider the last one, namely, the case $i = k-1$, it is easy to know that

$$J_k \ll \tilde{H}_k \tilde{Z}_k P S_{s-1}(P_{k-1}),$$

So, it has

$$P(\tilde{H}_{k-1} \tilde{Z}_{k-1})^2 Z_k S_{s-1}(P_{k-1}) = \tilde{H}_{k-1} \tilde{Z}_{k-1} \cdot \tilde{H}_k \tilde{Z}_k P S_{s-1}(P_{k-1})$$

And,

$$H_k = 1, \quad \text{i.e. } \theta_k = 1/k. \tag{3.7}$$

In general, it has

$$U_{j-1} = V_{j-1} \Rightarrow P(\tilde{H}_{j-1} \tilde{Z}_{j-1})^2 Z_j S_{s-1}(P_j) = \tilde{H}_{j-1} \tilde{Z}_{j-1} J_j$$
$$\Rightarrow P(\tilde{H}_{j-1} \tilde{Z}_{j-1}) Z_j S_{s-1}(P_j) = J_j (\approx U_j)$$
$$= S_{s-1}(P_j)^{1/2} Z_{j+1}^{(2s-3)/2} \left( P(\tilde{H}_j \tilde{Z}_j)^2 Z_{j+1} S_{s-1}(P_{j+1}) \right)^{1/2}$$
$$\Rightarrow P S_{s-1}(P_j) = Z_{j+1}^{(2s-2)} (H_j)^2 S_{s-1}(P_{j+1})$$
$$\Rightarrow Z_j^{2k} = S_{s-1}(P_{j+1}) S_{s-1}(P_j)^{-1} P Z_{j+1}^{(2s-2)} = P_{j+1}^{\lambda_{s-1}} P_j^{-\lambda_{s-1}} P Z_{j+1}^{(2s-2)} = Z_{j+1}^{-\lambda_{s-1}} P Z_{j+1}^{(2s-2)}$$
$$\Rightarrow P^{2k\theta_j} = P^{((2s-2)-\lambda_{s-1})\theta_{j+1}+1}.$$

i.e.

$$\theta_j = a \cdot \theta_{j+1} + b, \quad a = \left( (2s-2) - \lambda_{s-1} \right)/2k, \ b = 1/2k. \tag{3.8}$$

It may be shown directly that (3.8) is also true for $j = 1$.

As in [8], define $\Delta(s) = \lambda_s - (2s-k)$, from the recurrence (3.8), and with initial value (3.7), it follows

$$\theta_j = \frac{1}{k+\Delta(s-1)} + \left(\frac{1}{k} - \frac{1}{k+\Delta(s-1)}\right)\left(\frac{k-\Delta(s-1)}{2k}\right)^{k-j}. \tag{3.9}$$

Especially,

$$\theta = \theta_1 = \frac{1}{k+\Delta(s-1)} + \left(\frac{1}{k} - \frac{1}{k+\Delta(s-1)}\right)\left(\frac{k-\Delta(s-1)}{2k}\right)^{k-1}. \tag{3.10}$$

Then combine the result of section 2, here $P_0 = \tilde{P}$, it has

$$\tilde{P}^{\lambda_s} \ll Z^{2s-2}\sum_q (\sum_p T_{p,q}) \ll Z^{2s-1} \cdot ZPS_{s-1}(P) \ll Z^{2s} \cdot P^{1+\lambda_{s-1}}.$$

i.e.

$$\lambda_s \le \frac{\lambda_{s-1}}{1+\theta} + \frac{(1+2s\theta)}{1+\theta}. \tag{3.11}$$

The recursive inequality above resemble (2.6), but here $\theta$ is variable as $s$. (3.11) may be written as another form

$$\Delta(s) = \frac{\Delta(s-1)}{1+\theta} + \frac{(k\theta-1)}{1+\theta}. \tag{3.12}$$

Although it is not easy to deduce the exact expression of $\Delta(s)$ from (3.12) and (3.10), we may give an approximation

$$\Delta(s) \ll k \cdot \exp\left(-\frac{2(s-1)}{k+1}\right). \tag{3.13}$$

For the second part of $\theta$ in (3.10) will be much less than the first part, as $k$ is larger, so it may has

$$\theta \doteq \frac{1}{k+\Delta(s-1)}. \tag{3.10'}$$

Substituting (3.10') into (3.12)

$$\Delta(s) \doteq \Delta(s-1) - \frac{2\Delta(s-1)}{(k+1+\Delta(s-1))} = \Delta(s-1)\left(1 - \frac{2}{(k+1+\Delta(s-1))}\right)$$

Or,

$$\frac{\Delta(s)}{\Delta(s-1)} \doteq 1 - \frac{2}{(k+1)}\frac{1}{(1+\Delta(s-1)/(k+1))} \tag{3.14}$$

Let

$$\Delta(s) = (k-1)\prod_{2\le i\le s}\left(1-\frac{2\beta_i}{k+1}\right)$$

From (3.14), it has

$$\left(1-\frac{2\beta_s}{k+1}\right) \le 1-\frac{2}{(k+1)}\left(1+\prod_{2\le i\le s-1}\left(1-\frac{2\beta_i}{k+1}\right)\right)^{-1},$$

and

$$1-\beta_s \le \prod_{i\le s-1}\left(1-\frac{2\beta_i}{k+1}\right) \le \exp\left(-\sum_{i\le s-1}\left(\frac{2\beta_i}{k+1}\right)\right).$$

i.e.

$$\log(1-\beta_s) \le -\sum_{i\le s-1}\left(\frac{2\beta_i}{k+1}\right). \tag{3.15}$$

We may take the solution of the following differential equation as approximate solution of (3.15),

$$\frac{1}{(1-y)}\frac{dy}{dx} = \frac{2}{k+1}y.$$

It follows that $y(x) = 1-(e^{2(x-2)/(k+1)}+1)^{-1}$, and

$$\sum_{2\le i\le s}\beta_i = \sum_{2\le x\le s} y(x) > s-1-\frac{k+1}{2}\log\left(1+e^{2/(k+1)}\right).$$

And so

$$\Delta(s) = (k-1)\prod_{2\le i\le s}\left(1-\frac{2\beta_i}{k+1}\right)$$

$$\le (k-1)\exp\left(\sum_{2\le i\le s}-\frac{2\beta_i}{k+1}\right)$$

$$\le 2k\exp\left(-\frac{2(s-1)}{k+1}\right). \tag{3.16}$$

From section 2, we know

$$G(k) \le 3+2u+2\left[\frac{\Delta(u)}{2\hat{\sigma}}\right].$$

With (3.16) and take $u = 1+\left\lceil\frac{k+1}{2}\log\left(1/\hat{\sigma}\right)\right\rceil$, (1.2) is followed.